\documentclass[12pt,english]{article}
\usepackage[margin=0.75in]{geometry}
\usepackage[utf8]{inputenc}
\usepackage[english]{babel}
\usepackage{amsthm}
\usepackage{graphicx}
\usepackage{url}
\usepackage{hyperref}
\usepackage{amsmath}
\usepackage{listings}
\usepackage{amsfonts}
\usepackage{xcolor}
\usepackage{float}
\newtheorem{theorem}{Theorem}

\def\XXint#1#2#3{{\setbox0=\hbox{$#1{#2#3}{\int}$}
     \vcenter{\hbox{$#2#3$}}\kern-.5\wd0}}

\begin{document}

\title{Analytical continuation of Euler prime product for $\Re(s)>\tfrac{1}{2}$ assuming (RH)}

\author{Artur Kawalec}

\date{}
\maketitle

\begin{abstract}
We analytically continue the Euler prime product for $\Re(s)>\tfrac{1}{2}$ (except for its pole $s=1$) assuming (RH) by introducing a new factor to the Euler product. We also discuss how to recover the Mertens's 3rd Theorem at $s=1$ case, and how to apply the same technique to analytically continue other similar Euler products. In the last part, we also construct a simple script in Pari/GP to compute the Euler product and verify the calculations numerically.
\end{abstract}

\section{Introduction}
Let $p=\{2,3,5,7,\ldots\}$ be a sequence of primes, the Euler product over primes defines the Riemann zeta function

\begin{equation}\label{eq:1}
\zeta(s)=\prod_{p}\left(1-\frac{1}{p^s}\right)^{-1}
\end{equation}
which is absolutely convergent for $\Re(s)>1$, where $s=\sigma+it$ is a complex variable. The Euler product shows in one of the simplest terms, how primes enter the zeta function, and by the virtue of its absolute convergence, it proves that the zeta function has no zeros for $\Re(s)>1$. At $s=1$, the Euler product diverges and the Mertens's 3rd theorem

\begin{equation}\label{eq:1}
\prod_{p\leq x}\left(1-\frac{1}{p}\right)^{-1}\sim e^{\gamma}\log x
\end{equation}
describes the asymptotic growth of the partial product terms at $s=1$, which grow as $\log x$ and $\gamma=0.57721566\ldots$ is the Euler's constant. There are interesting papers by Gonek about a hybrid Euler product [3][4]. And in this article, we analytically continue the Euler product for $\Re(s)>\tfrac{1}{2}$ domain ($s\neq 1$), as shown by next Theorem:

\begin{theorem}

\begin{equation}\label{eq:1}
\zeta(s)=\lim_{x\to \infty}e^{\operatorname{E_1}[(s-1)\log x]}\prod_{p\leq x}\left(1-\frac{1}{p^s}\right)^{-1}
\end{equation}
valid for $\Re(s)>\tfrac{1}{2}$ and $(s\neq 1)$ assuming (RH).
\end{theorem}
\noindent
where the exponential integral is defined by

\begin{equation}\label{eq:1}
\operatorname{E_1}(z)=\int_{z}^{\infty}\frac{e^{-t}}{t}dt, \quad |\arg(z)|<\pi
\end{equation}
valid for complex $z$ with a branch cut on $(-\infty, 0]$.

\section{Proof of Theorem 1}

To do this, we first uncover what prevents the analytical continuation of Euler product in the first place. If we consider taking the logarithm of Euler product

\begin{equation}\label{eq:1}
\log\zeta(s)=\sum_{n=1}^{\infty}\sum_{p}\frac{1}{np^{ns}}=\sum_{n=1}^{\infty}\frac{P(ns)}{n}
\end{equation}
where
\begin{equation}\label{eq:1}
P(s)=\sum_{p}\frac{1}{p^s}
\end{equation}
is defined as the prime zeta function, which is also absolutely convergent for $\Re(s)>1$. As a result, rearranging the series such that

\begin{equation}\label{eq:1}
\log\zeta(s)=P(s)+\sum_{n=2}^{\infty}\frac{P(ns)}{n}
\end{equation}
shows that $P(s)$ converges for $\Re(s)>1$, while the second series
\begin{equation}\label{eq:1}
\sum_{n=2}^{\infty}\frac{P(ns)}{n}
\end{equation}
converges for $\Re(s)>\tfrac{1}{2}$; therefore it is the $P(s)$ component that actually limits further convergence of Euler product to $\frac{1}{2}<\Re(s) < 1$, which is a region to the left side of the line $\Re(s)=1$. Technically, the Euler product (1) still holds on the line $\Re(s)=1$ when $(t \neq 0)$ as shown in [7, p.65]. In our previous article [5], we showed that

\begin{equation}\label{eq:1}
P(s)=\sum_{p\leq x}\frac{1}{p^s}+\operatorname{E_1}[(s-1)\log(x)]+O\left(x^{\frac{1}{2}-s}\log x\right)
\end{equation}
converges for $\Re(s)>\tfrac{1}{2}$ as $x\to \infty$ having a branch cut on $(\tfrac{1}{2},1]$ assuming (RH). As a result, considering the partial Euler product

\begin{equation}\label{eq:1}
\log \left[\prod_{p\leq x}\left(1-\frac{1}{p^s}\right)^{-1}\right]=\sum_{p\leq x}\frac{1}{p^s}+\sum_{n=2}^{\infty}\frac{1}{n}\left(\sum_{p\leq x}\frac{1}{p^{ns}}\right)
\end{equation}
and then adding the exponential integral factor

\begin{equation}\label{eq:1}
\sum_{p\leq x}\frac{1}{p^s}+\sum_{n=2}^{\infty}\frac{1}{n}\left(\sum_{p\leq x}\frac{1}{p^{ns}}\right)+\operatorname{E_1}[(s-1)\log(x)]+O\left(x^{\frac{1}{2}-s}\log x\right)
\end{equation}
will extend the $P(s)$ component to $\Re(s)>\tfrac{1}{2}$ with a branch cut on $(\tfrac{1}{2},1]$ assuming (RH). And finally, when exponentiating back we will recover the Euler product (3) as

\begin{equation}\label{eq:1}
\zeta(s)=e^{\operatorname{E_1}[(s-1)\log x]}\prod_{p\leq x}\left(1-\frac{1}{p^s}\right)^{-1}\left(1+O\Big(x^{\frac{1}{2}-s}\log x\Big)\right)
\end{equation}
and taking limit $x\to\infty$.

\section{A note on Euler product at $s=1$}

At $s=1$, we can recover the Mertens's 3rd Theorem. If we consider the Riemann zeta near the pole as $s\to 1$, then we expect the limit

\begin{equation}\label{eq:1}
\lim_{s\to 1^{+}}\zeta(s)(s-1)=1
\end{equation}
Now substituting the Euler product we must have
\begin{equation}\label{eq:1}
\lim_{s\to 1^{+}}(s-1)e^{\operatorname{E_1}[(s-1)\log x]}\prod_{p\leq x}\left(1-\frac{1}{p^s}\right)^{-1}= 1
\end{equation}
The series expansion of exponential integral is

\begin{equation}\label{eq:1}
\operatorname{E_1}(z)=-\gamma-\log(z)-\sum_{k=1}^{\infty}\frac{(-z)^k}{k\cdot k!}
\end{equation}
and so, near $s\to 1$, we have
\begin{equation}\label{eq:1}
\operatorname{E_1}[(s-1)\log x]=-\gamma-\log(s-1)-\log\log x+o(1)
\end{equation}
therefore

\begin{equation}\label{eq:1}
\frac{e^{-\gamma}}{\log x}\prod_{p\leq x}\left(1-\frac{1}{p}\right)^{-1}=1
\end{equation}
and 
\begin{equation}\label{eq:1}
\prod_{p\leq x}\left(1-\frac{1}{p}\right)^{-1}\sim e^{\gamma}\log x.
\end{equation}

\section{A note on other related Euler Products}

The first example we investigate is the reciprocal zeta function
\begin{equation}\label{eq:17}
\frac{1}{\zeta(s)}=\sum_{n=1}^{\infty}\frac{\mu(n)}{n^s}=\prod_{p}\left(1-\frac{1}{p^s}\right)
\end{equation}
The series involving the M\"obius function converges for $\Re(s)>\frac{1}{2}$ assuming (RH), while the Euler product is valid for $\Re(s)>1$. Therefore, simply adding the exponential correction factor

\begin{equation}\label{eq:1}
\frac{1}{\zeta(s)}=\lim_{x\to \infty}e^{-\operatorname{E_1}[(s-1)\log x]}\prod_{p\leq x}\left(1-\frac{1}{p^s}\right)
\end{equation}
extends it for $\Re(s)>\frac{1}{2}$ assuming (RH), where it has a zero at $s=1$. 

And in another example, a closely related Euler product

\begin{equation}\label{eq:17}
\frac{\zeta(2s)}{\zeta(s)}=\prod_{p}\left(1+\frac{1}{p^{s}}\right)^{-1}
\end{equation}
also is valid for $\Re(s)>1$. And the log expansion gives
\begin{equation}\label{eq:17}
\log\left[\frac{\zeta(2s)}{\zeta(s)}\right]=-P(s)+\sum_{n=2}^{\infty}(-1)^n  \frac{P(ns)}{n}
\end{equation}
which shows that it converges for $\Re(s)>1$ because of being limited $P(s)$ component (6) valid for $\Re(s)>1$, so if we extend that as before, then we get

\begin{equation}\label{eq:17}
\frac{\zeta(2s)}{\zeta(s)}=e^{-\operatorname{E_1}[(s-1)\log x]}\prod_{p\leq x}\left(1+\frac{1}{p^{s}}\right)^{-1}
\end{equation}
gives analytical continuation to $\Re(s)>\frac{1}{2}$.

\section{Numerical computation}

In this Section, we will numerically verify this Euler product. In Listing $1$, we write a script to calculate the Euler product in Pari/GP software package [6], where we utilize the built-in exponential integral function $E_1(Z)$ as $\textbf{eint1(z)}$. And in Fig. 1, we plot the zeta function (blue trace) and compare with the analytically continued Euler product equation (3) in the range $\tfrac{1}{2}<\Re(s)<2$ (dotted red trace) for $x=10^6$. The result is that we see a perfect match, but near $s=0.5$ the Euler product starts deviating a little bit. And in Fig. 2 we plot the modulus of the analytically continued Euler product on a vertical line with real part $\sigma=0.8$ and $t=0$ to $t=50$ at limit variable $x=10^3$, where we also see a near perfect match with the zeta function. In Fig. 3 we make the same plot but with $\sigma=0.55$, which is near the critical line, and we start seeing a few oscillations that show up. So in Fig. 4 we re-plot the $\sigma=0.55$ case again, but increase $x=10^5$, and we observe that these oscillations smooth out even more. Eventually as $x\to\infty$  these oscillations will smooth out and will approach the zeta function.

\lstset{language=C,deletekeywords={for,double,return},caption={A Pari/Gp code excerpt for computing Euler product equation (3)},label=DescriptiveLabel,captionpos=b,showstringspaces=false}
\begin{lstlisting}[frame=single]
\\ Set Limit variable x
x = 10^3;

\\ Define Euler Prime Product by equation (3)
EulerProd(s) = exp(eint1((s-1)*log(x)))*
    prod(n=1,primepi(x),(1-1.0/prime(n)^s)^(-1));

\\ Fig.1 Plot
ploth(s=0.501, 2, [zeta(s), real(EulerProd(s))])

\\ Fig.2-4 Plots
a = 0.8;
ploth(t=0, 50, [abs(zeta(a+I*t)), abs(EulerProd(a+I*t))])

\end{lstlisting}

\texttt{Email: art.kawalec@gmail.com}

\renewcommand{\figurename}{Figure}
\begin{figure}[H]
  \renewcommand{\thefigure}{1}%
  \centering
  \includegraphics[width=175mm]{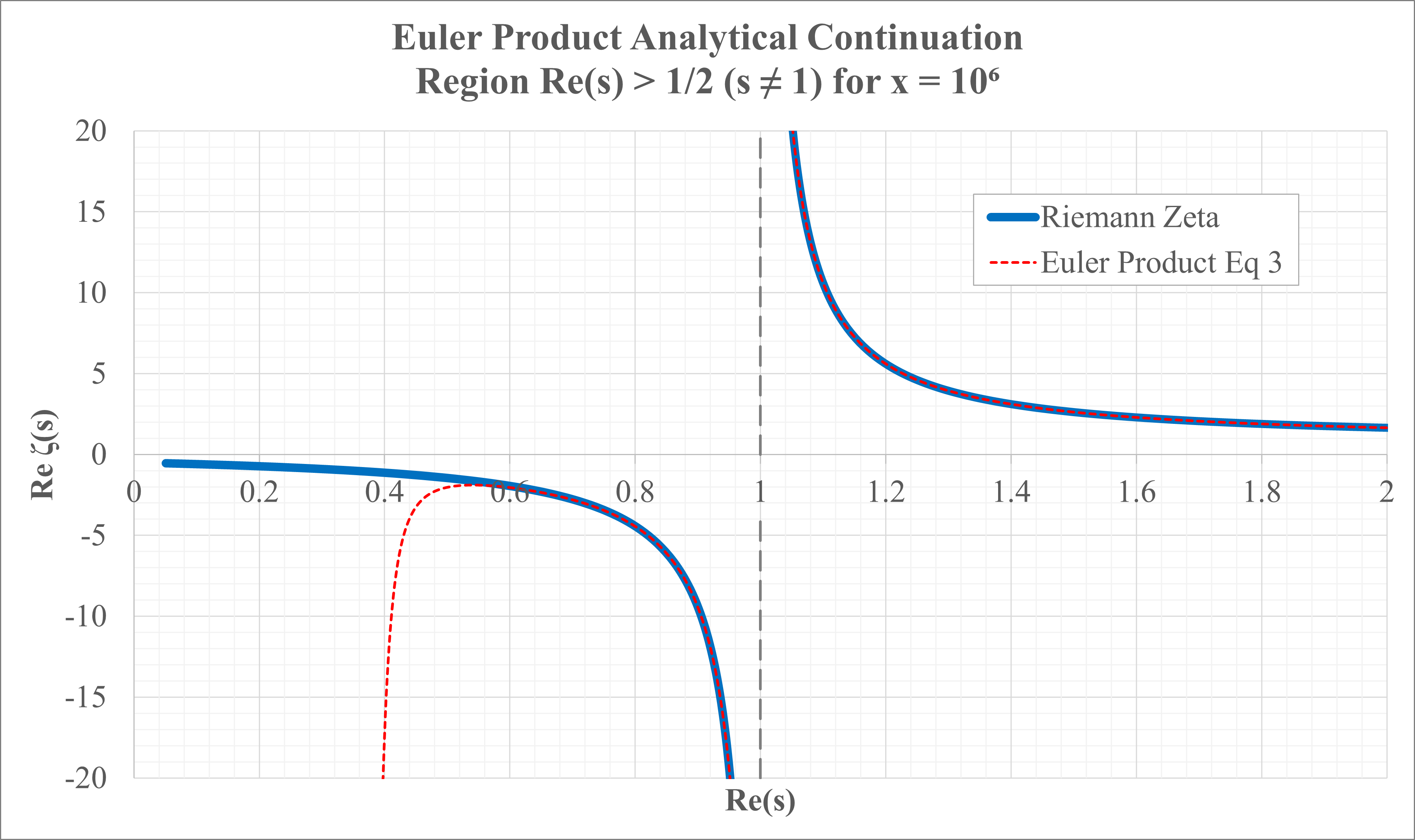}\\
  \caption{A plot of the analytically continued Euler product for real $s$ for $x=10^6$}\label{1}
\end{figure}

\renewcommand{\figurename}{Figure}
\begin{figure}[H]
  \renewcommand{\thefigure}{2}%
  \centering
  \includegraphics[width=175mm]{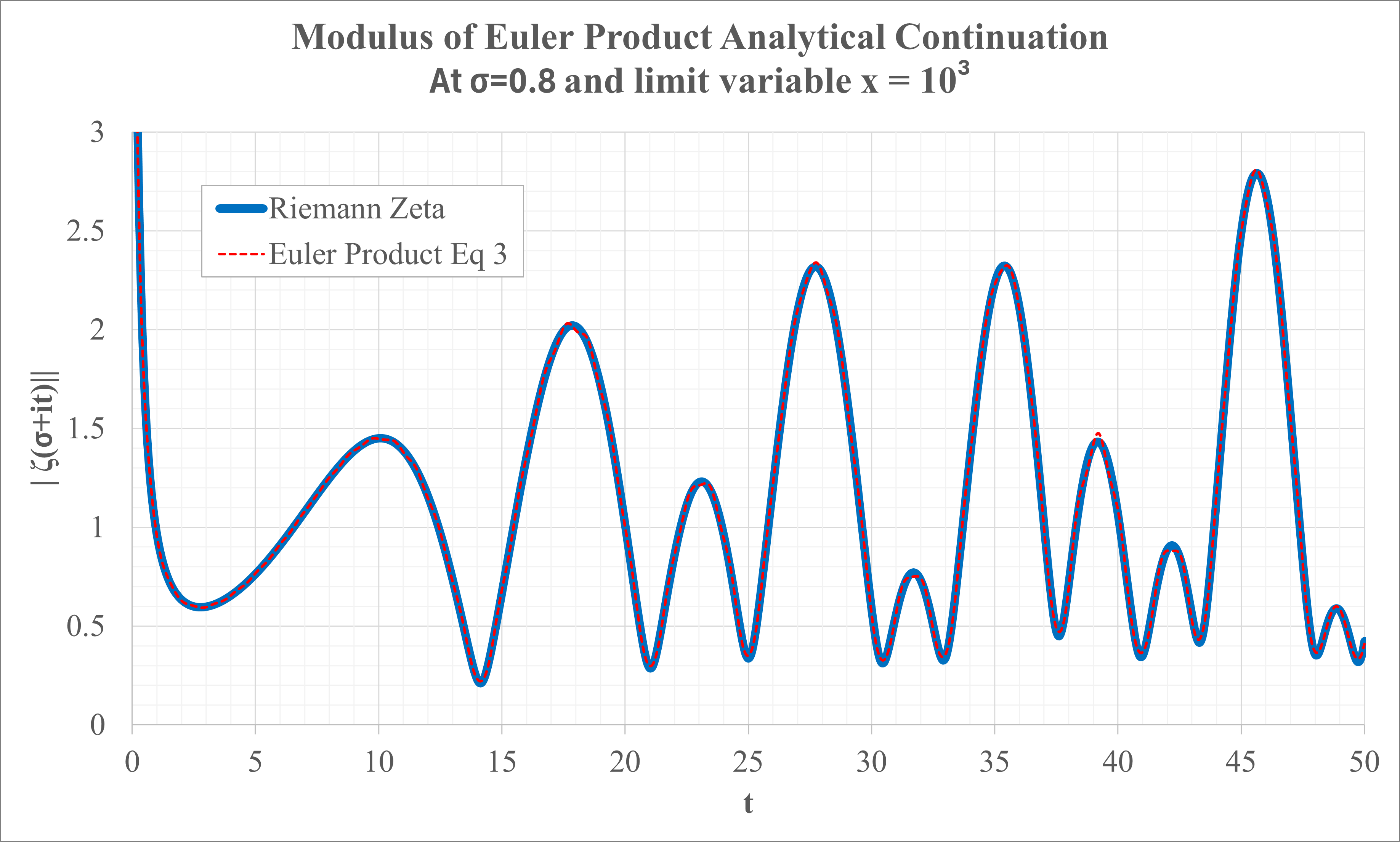}\\
  \caption{A plot of the analytically continued Euler product at $\sigma=0.8$ for $x=10^3$}\label{1}
\end{figure}

\renewcommand{\figurename}{Figure}
\begin{figure}[H]
  \renewcommand{\thefigure}{3}%
  \centering
  \includegraphics[width=175mm]{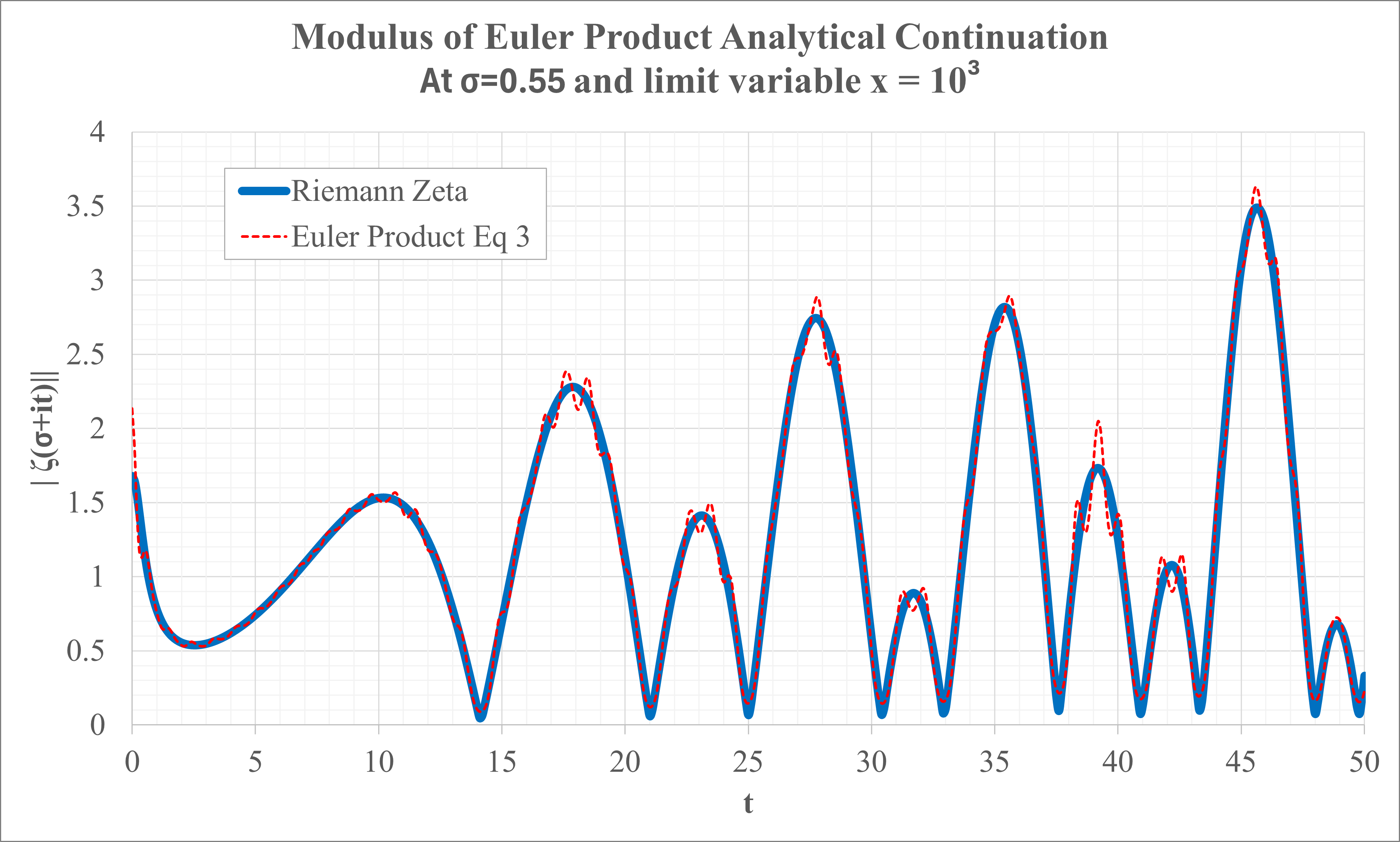}\\
  \caption{A plot of the analytically continued Euler product at $\sigma=0.55$ for $x=10^3$}\label{1}
\end{figure}

\renewcommand{\figurename}{Figure}
\begin{figure}[H]
  \renewcommand{\thefigure}{4}%
  \centering
  \includegraphics[width=175mm]{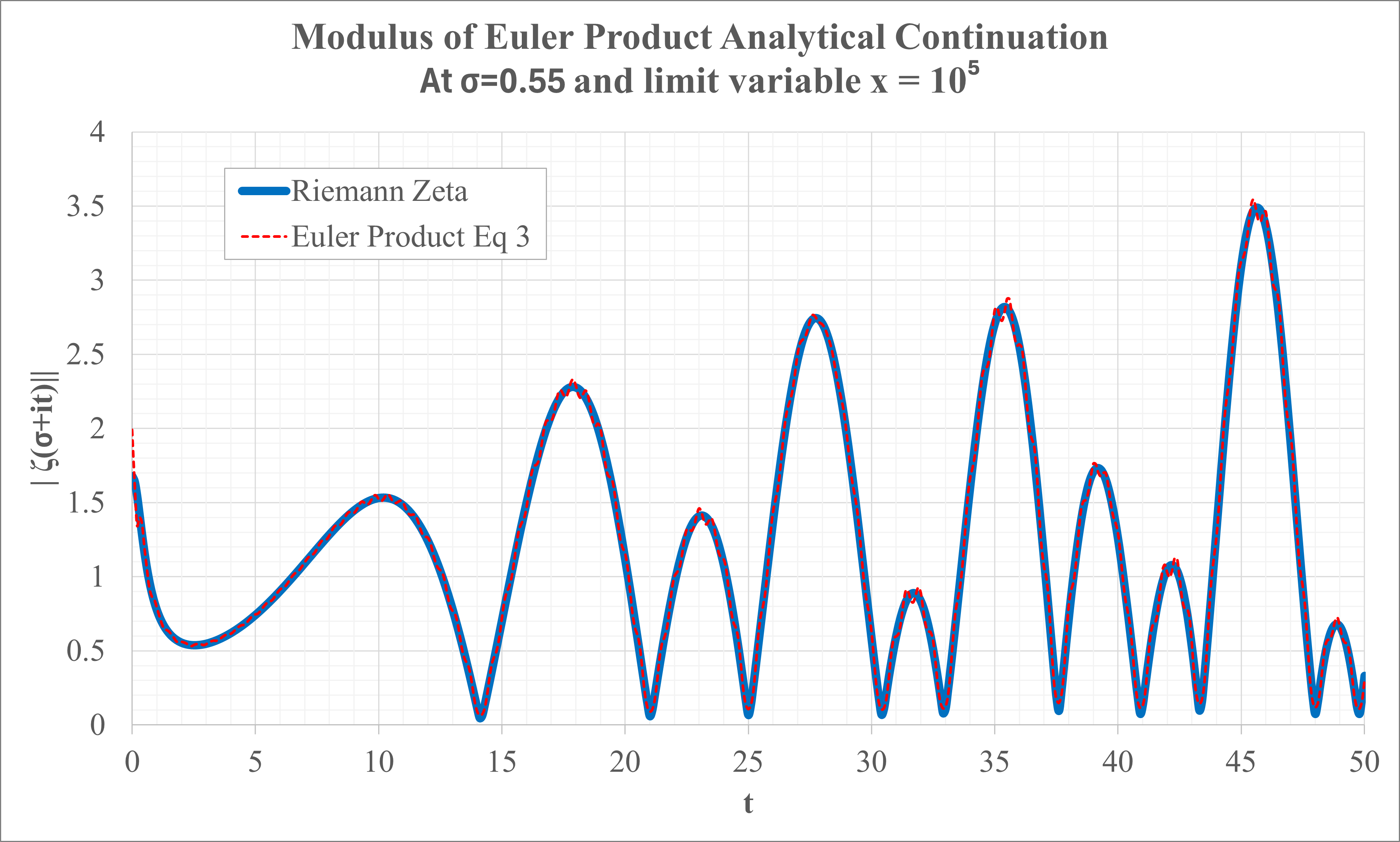}\\
  \caption{A plot of the analytically continued Euler product at $\sigma=0.55$ for $x=10^5$}\label{1}
\end{figure}

\end{document}